# FILOLOGIA CLASSICA
Vittorio Pizzigoni

Nel margine della sua copia dell'*Arithmetica* di Diofanto Fermat scrive: «È impossibile scrivere un cubo come somma di due cubi o una quarta potenza come somma di due quarte potenze o, in generale, nessun numero che sia una potenza maggiore di due può essere scritto come somma di due potenze dello stesso valore»[1]. Parafrasando le parole di Fermat affermano che l'equazione $X^n + Y^n = Z^n$ non ammette soluzioni intere per *n*=3, non ammette soluzioni intere per *n*=4, e che «generalizzando» non ammette soluzioni per $n > 4$.

A margine di questo enunciato aggiunge la frase divenuta famosa: «Dispongo di una meravigliosa dimostrazione di questo teorema che non può essere contenuta nel margine troppo stretto della pagina»[2].

Cerchiamo di penetrare il senso di queste parole con gli strumenti filologici. Innanzitutto quale era lo stato dell'arte al momento in cui Fermat scriveva questa nota? È sicuro che Fermat dimostrò che l'equazione non ammette soluzioni per *n*=4, e sembra molto probabile che dimostrò anche quello in cui *n*=3. Queste due dimostrazioni, poi sviluppate nel dettaglio da Eulero, presentano alcune differenze ma risultano concettualmente simili: entrambe utilizzano una dimostrazione per assurdo ricorsiva, che si può indicare come «metodo della discesa infinita». Si suppone che l'equazione data, sia essa $X^3 + Y^3 = Z^3$ oppure $X^4 + Y^4 = Z^4$, ammetta una soluzione intera e si dimostra che questo implica l'esistenza di una soluzione più piccola formalmente simile a quella iniziale; è quindi possibile ripetere il procedimento arrivando a trovare una soluzione ancora più piccola, e procedendo così all'infinito si dimostra l'assurdità dell'assunto iniziale.

Sorprende allora che Fermat nella sua enunciazione del problema usi il termine *«generaliter»* e se prestiamo la dovuta attenzione a questo termine dobbiamo ipotizzare che Fermat pensasse di poter 'generalizzare' a tutte le infinite potenze superiori a 4 ai due casi in cui *n*=3 e *n*=4. Partendo da questa considerazione cerchiamo di ricostruire il pensiero di Fermat.

Se immaginiamo di sviluppare sul piano le potenze superiori a 4 di un numero ci accorgiamo che esse presentano una forma ricorsiva: dato un numero qualsiasi ogni sua potenza pari è *formalmente* simile ad ogni altra sua potenza pari, così come ogni sua potenza dispari è *formalmente* simile ad ogni altra sua potenza dispari, per aiutarci possiamo osservare la seguente illustrazione che esemplifica le potenze del numero 3.

---

[1] «Cubem autem in duos cubos, aut quadratoquadratum in duos quadratoquadratos, et generaliter nullam in infinitum ultra quadratum potestatem in duos eiusdem nominis fas est dividere». Cfr. Simon Singh, *L'ultimo teorema di Fermat*, BUR Biblioteca Univ. Rizzoli, Milano 1999, p. 88.

[2] «Cuius rei demonstrationem mirabilem sane detexi hanc marginis exiguitas non caperet». Cfr. S. Singh, *L'ultimo teorema...*, cit., p. 89.

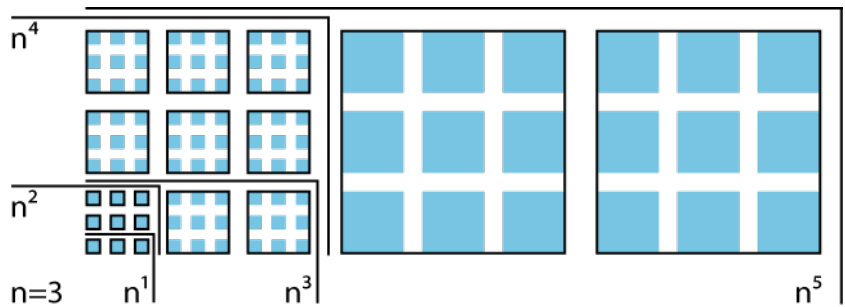

Sviluppo sul piano delle prime potenze del numero 3.

Grazie a tale similitudine formale è possibile immaginare di suddividere in due equazioni l'equazione $X^n + Y^n = Z^n$ :

1)  $A^3 + B^3 = C^3$  (per ogni $n$ dispari superiore a 3)
2)  $A^4 + B^4 = C^4$  (per ogni $n$ pari superiore a 4)

dove A, B, C sono '*sovraunità*' ricavate dalle potenze $X^n$, $Y^n$, $Z^n$.
Per esempio: se X=2, $n$=6, e quindi $2^6$=64 possiamo immaginare di scrivere 16*a* - dove ognuna delle 16 '*sovraunità*' *a* è formata da 4 unità - oppure $\{2a\}^4$. Più in generale possiamo scrivere ogni $X^n$, $Y^n$, $Z^n$ come $\{Xa\}^3$, $\{Yb\}^3$, $\{Zc\}^3$ nel caso $n$ sia dispari; o come $\{Xa\}^4$, $\{Yb\}^4$, $\{Zc\}^4$ nel caso $n$ sia pari.
Per brevità adottiamo la notazione di 1) e 2).
Data la coprimarietà di X, Y, Z e poiché le '*sovraunità*' A, B, C contengono sempre una potenza di X, Y, Z anche A, B e C sono coprimi. In questo modo possiamo ricondurre ogni equazione in cui $n$ sia superiore a 4 ai casi in cui $n$=3 e $n$=4. Poiché per $n$=3 e per $n$=4 le precedenti dimostrazioni portano ad una ricorsività discendente siamo sicuri che, per quanto grande sia il numero di unità contenute nelle '*sovraunità*' A, B, C, esso verrà sempre superato dalla infinita reiterazione permessa dalle prime due dimostrazioni. Pertanto si conclude che per ogni potenza superiore a due l'equazione $X^n + Y^n = Z^n$ non ammette soluzioni intere.
Questa dimostrazione pur potendo non trovare spazio nel margine dell'*Arithmetica* posseduta da Fermat, sembra abbastanza corta da avvicinarsi a quella brevità che egli stesso lascia intendere nella sua pungente annotazione.